\documentclass[11pt]{article}
\usepackage{psfig,epsfig,amssymb}
\usepackage{color}
\usepackage[colorlinks=true,
linkcolor=webgreen,
filecolor=webbrown,
citecolor=webgreen]{hyperref}

\definecolor{webgreen}{rgb}{0,.5,0}
\definecolor{webbrown}{rgb}{.6,0,0}

\def\slfrac#1#2{\hbox{\kern.1em %
 \raise.5ex\hbox{\the\scriptfont0 #1}\kern-.11em %
 /\kern-.15em\lower.25ex\hbox{\the\scriptfont0 #2}}}

\newcommand{\eeq}{\end{equation}}
\newcommand{\beql}[1]{\begin{equation}\label{#1}}

\makeatletter
\def\@sect#1#2#3#4#5#6[#7]#8{\ifnum #2>\c@secnumdepth
     \def\@svsec{}\else
     \refstepcounter{#1}\edef\@svsec{\csname the#1\endcsname.\hskip .75em }\fi
     \@tempskipa #5\relax
      \ifdim \@tempskipa>\z@
        \begingroup #6\relax
          \@hangfrom{\hskip #3\relax\@svsec}{\interlinepenalty \@M #8\par}%
        \endgroup
       \csname #1mark\endcsname{#7}\addcontentsline
         {toc}{#1}{\ifnum #2>\c@secnumdepth \else
                      \protect\numberline{\csname the#1\endcsname}\fi
                    #7}\else
        \def\@svsechd{#6\hskip #3\@svsec #8\csname #1mark\endcsname
                      {#7}\addcontentsline
                           {toc}{#1}{\ifnum #2>\c@secnumdepth \else
                             \protect\numberline{\csname the#1\endcsname}\fi
                       #7}}\fi
     \@xsect{#5}}
\def\@begintheorem#1#2{\it \trivlist \item[\hskip \labelsep{\bf #1\ #2.}]}
\makeatother

\begin{document}
\begin{center}
{\Large{\bf The On-Line Encyclopedia of Integer Sequences}} \\
\vspace{1\baselineskip}
{\em N. J. A. Sloane}\footnotemark \\
\vspace{.25\baselineskip}
\vspace{2\baselineskip}
\vspace{.5\baselineskip}
\end{center}
\footnotetext[1]{Neil J. A. Sloane is with AT\&T Shannon Labs, Florham Park, NJ.
His email address is
\href{mailto:njas@research.att.com}{njas@research.att.com}.}

\setlength{\baselineskip}{1.5\baselineskip}

\section{Introduction}
This article gives a brief introduction to the
On-Line Encyclopedia of Integer Sequences (or OEIS).
The OEIS is a database of nearly $90,000$ 
sequences of integers, arranged lexicographically.
The entry for a sequence lists the initial terms (50 to 100, if available),
a description, formulae,
programs to generate the sequence,
references, links to relevant web pages, and other information.

\section{To consult the database}
Since 1996, an electronic version \cite{OEIS}
has been accessible via the Internet, at the URL
$$
\htmladdnormallink{\mbox{http://www.research.att.com/$\sim$njas/sequences/}}
{http://www.research.att.com/~njas/sequences/}.
$$
If a list of numbers is entered there,
the reply will display the entries for all matching sequences.

For example, suppose you were trying to count the ways to insert parentheses into a string of $n$ letters so that the parentheses are balanced and there are at least two letters inside each pair of parentheses.
The outer pair of parentheses is to be ignored.
For $n=3$ and 4 there are respectively 3 and 11 possibilities:
$$
\begin{array}{ll}
n=3: & abc, (ab)c , a(bc); \\
n=4: & abcd, (ab)cd, a(bc)d, ab(cd), (ab)(cd), (abc)d, \\
& a(bcd), ((ab)c)d, (a(bc))d, a((bc)d), a(b(cd)) \,.
\end{array}
$$
Further work shows that for $n=1, \ldots, 5$ the numbers
are 1, 1, 3, 11, 45.
Entering these into the web page produces nine matching sequences,
but they are sorted, with the most probable match appearing first.
Indeed, this entry tells you that these are the numbers
(sequence
\htmladdnormallink{A1003}{http://www.research.att.com/cgi-bin/access.cgi/as/njas/sequences/eisA.cgi?Anum=A001003}) arising from ``Schr\"{o}der's second problem'', and are also known as ``super-Catalan numbers''.

\begin{figure}[htb]
\caption{Part of the reply when the sequence
1, 1, 3, 11, 45 is submitted to the On-Line Encyclopedia.
Many references, links and comments
have been omitted to save space.}

\vspace*{+.25in}
\centerline{\psfig{file=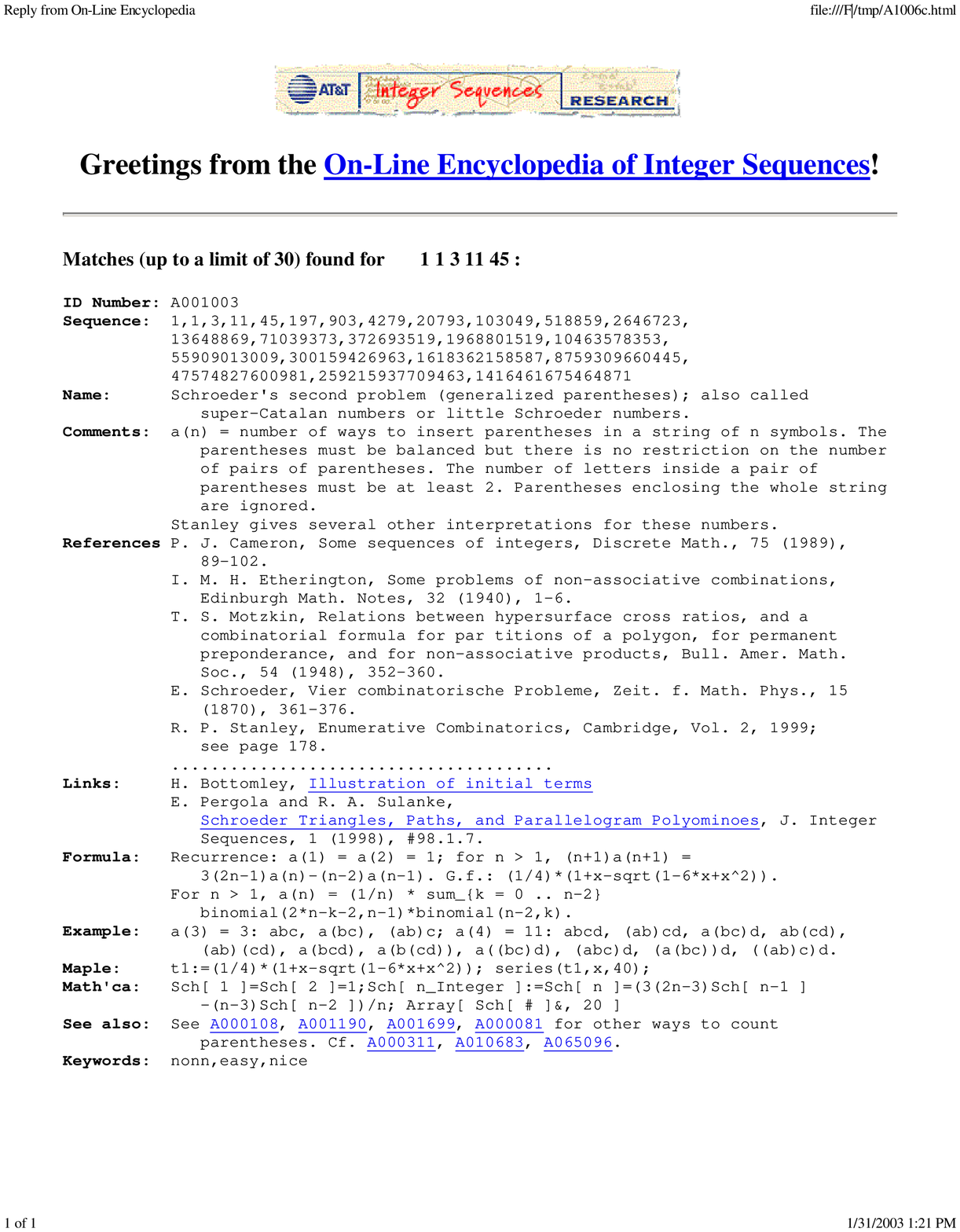,width=6in}}

\label{F1}
\end{figure}

The reply (an abridged version is shown in Figure \ref{F1}) gives 21 references, ranging from Schr\"{o}der (1870) \cite{Sch70} to articles
published electronically in the last few years.
There is an explicit formula:
$$a(n) = \frac{1}{n}
\sum_{k=0}^{n-2} \left( \begin{array}{c}
2n-k-2 \\ n-1
\end{array}\right)
\left( \begin{array}{c}
n-2 \\ k
\end{array}\right)~ , ~~ n > 1 \,,
$$
a recurrence:
$$a(1) = a(2) =1 \,,$$
$$(n+1) a(n+1) = 3 (2n-1) a(n) - (n-2) a(n-1) ~ , ~~ n > 1 \,,
$$
programs to produce the sequence in Maple and Mathematica, and much more.

There's no other reference work that will carry out this kind of search.

The Encyclopedia can also be consulted via email.
There are two addresses.
Sending email to 
\href{mailto:sequences@research.att.com}{sequences@research.att.com}
with a line in the body of the message saying
$$\mbox{lookup~~1~~1~~3~~11~~45}$$
will trigger the same search that the web page performs,
only now the results are sent, almost immediately, via email.
Superseeker
(\href{mailto:superseeker@research.att.com}{superseeker@research.att.com})
carries out a more sophisticated analysis and 
tries hard to find an explanation for the sequence, even if it is not in the database.
If the simple lookup fails, Superseeker carries out many other tests, including:
\begin{itemize}
\item
applying over 130 transformations to the sequence, including the binomial, Euler, M\"{o}bius, etc. transforms \cite{BeSl95},
and checking to see if the result is in the database
\item
applying Pad\'{e} approximation methods to try, for example, to 
express the $n^{\rm th}$ term
as a rational function of $n$ (using the ``gfun'' package of Salvy and Zimmermann \cite{SaZi94}, the ``guesss'' program 
of Derksen \cite{Der94} and the ``RATE'' program of Kratthentaler \cite{Kra00})
\item
checking to see if changing one or two characters produces a sequence in the database.
\end{itemize}

\noindent
Since Superseeker carries out a nontrivial amount of calculation,
users are asked to submit only one sequence per hour.

The electronic version of the Encyclopedia had its origins in the books \cite{HIS} (1973) and \cite{EIS} (1995).
Disk space is cheap, and the present incarnation (excluding illustrations) contains 
about $72$ times as much data as the 1995 book.
The history of the Encyclopedia is described in more detail in \cite{Slo98}.

\section{Applications}
Most people use the Encyclopedia to identify a sequence, as illustrated above.
It has been around long enough so that there is a good chance that your sequence will be there.
If not, you will see a message encouraging you to submit it.

Most of these applications are unspectacular,
akin to looking up a word in a dictionary (cf. \cite{Cip94}).  One encounters
a sequence in the middle of a calculation, perhaps
$$1~~2~~4~~6~~10~~12~~16~~18~~22~~28~~30~~\ldots ~,$$
and one wants to know quickly what it is --- preferably
a formula for the $n$-th term
(in this case it is probably prime$(n) - 1$, \htmladdnormallink{A6093}{http://www.research.att.com/cgi-bin/access.cgi/as/njas/sequences/eisA.cgi?Anum=A006093})
or a recurrence.  Successful applications of this type
usually go unremarked. Some are more dramatic: there is 
a web page\footnote{
http://www.research.att.com/$\sim$njas/sequences/cite.html
}
that lists several hundred articles that acknowledge
help from the OEIS.
One quotation will serve to illustrate this.
Emeric Deutsch of Poytechnic University, Brooklyn,
said in a recent email message: ``... your database is
invaluable. For example, for a certain sequence $a_n$,
using Maple I found the first 100 or so
indices $i$ for which $a_i$ is odd. Only the OEIS could tell
me that the sequence of these $i$'s is a known sequence
related to the Thue-Morse sequence. Of course, this had
to be followed by further reading and proof.''

The other main application is to find out the current status of 
work on a problem --- for example,
the search for Mersenne primes (see \htmladdnormallink{A43}{http://www.research.att.com/cgi-bin/access.cgi/as/njas/sequences/eisA.cgi?Anum=A000043}), the
enumeration of Hadamard matrices (\htmladdnormallink{A7299}{http://www.research.att.com/cgi-bin/access.cgi/as/njas/sequences/eisA.cgi?Anum=A007299}), Latin squares (\htmladdnormallink{A315}{http://www.research.att.com/cgi-bin/access.cgi/as/njas/sequences/eisA.cgi?Anum=A000315}) or meanders (\htmladdnormallink{A5316}{http://www.research.att.com/cgi-bin/access.cgi/as/njas/sequences/eisA.cgi?Anum=A005316}) of order $n$, the
latest information about the decimal expansion of $\pi$ (\htmladdnormallink{A796}{http://www.research.att.com/cgi-bin/access.cgi/as/njas/sequences/eisA.cgi?Anum=A000796}) or,
better, its continued fraction expansion (\htmladdnormallink{A1203}{http://www.research.att.com/cgi-bin/access.cgi/as/njas/sequences/eisA.cgi?Anum=A001203}).

Of course people trying to solve puzzles or IQ tests find the database useful.
\htmladdnormallink{A5228}{http://www.research.att.com/cgi-bin/access.cgi/as/njas/sequences/eisA.cgi?Anum=A005228} is a classic:
$$1~~3~~7~~12~~18~~26~~35~~45~~56~~69~~83~~\ldots ~ $$

There are also some less obvious applications.
One is in simplifying complicated integer-valued expressions.
You might, for example, have encountered the sum
$$\sum_{k=0}^n \left( \begin{array}{c}
4n+1 \\ 2n-2k \end{array}\right)
\left( \begin{array}{c}
n+k \\ k \end{array}\right) \,.
$$
There are powerful methods for evaluating such sums
\cite{NePe97}, \cite{PWZ96},
but it doesn't take long to work out the first few terms:
1, 12, 240, 5376, and to look them up in the database.
In this case you would have been lucky.
The reply suggests that this is sequence \htmladdnormallink{A6588}{http://www.research.att.com/cgi-bin/access.cgi/as/njas/sequences/eisA.cgi?Anum=A006588},
$4^n \left( \begin{array}{c} 3n \\n \end{array}\right)$, and supplies, with references, the binomial coefficient
identity you were hoping for.

Another application is in proving inequalities.
You might suspect that $\sigma (n) < n \sqrt{n}$ for $n>2$, where $\sigma$ is the
sum-of-divisors function (\htmladdnormallink{A203}{http://www.research.att.com/cgi-bin/access.cgi/as/njas/sequences/eisA.cgi?Anum=A000203}).
If the initial terms of $[n \sqrt{n} ] - \sigma (n)$
(where $[\,]$ denotes the ``floor'' function)
are submitted to the database, the reply suggests that this is \htmladdnormallink{A55682}{http://www.research.att.com/cgi-bin/access.cgi/as/njas/sequences/eisA.cgi?Anum=A055682}, and points you to
a reference that gives a proof of your inequality.

I cannot resist mentioning sequence \htmladdnormallink{A57641}{http://www.research.att.com/cgi-bin/access.cgi/as/njas/sequences/eisA.cgi?Anum=A057641}, which gives the values of
$$\left[ H_n + \exp (H_n) \log (H_n ) \right] - \sigma (n)$$
for $n \ge 1$, where $H_n$ is the harmonic number
$\sum_{i=1}^n 1/i$.
Lagarias \cite{Lag02},
extending earlier work of Robin \cite{Rob84}, has shown that this sequence is nonnegative if and only if the Riemann hypothesis holds!

Although the database contains a number of sequences of
both of the above types, I have not made a systematic search through 
reference works such as \cite{Gou72}, \cite{HNT}
and it would be nice to get many more examples.

The database can also be used to save space when referring to particular sequences.
When introducing the Motzkin numbers,
for example, instead of giving the definition,
references and the first few terms it is simpler just to say 
``$\ldots$ the Motzkin numbers $M_n$ (sequence \htmladdnormallink{A1006}{http://www.research.att.com/cgi-bin/access.cgi/as/njas/sequences/eisA.cgi?Anum=A001006} of \cite{OEIS})''.

One can also search the database for sequences that mention a particular
name (Riemann, say),
and there is a separate alphabetical index,
useful for keeping track of all sequences on a certain topic
--- e.g. the entry for {\em groups}
lists abelian (\htmladdnormallink{A688}{http://www.research.att.com/cgi-bin/access.cgi/as/njas/sequences/eisA.cgi?Anum=A000688}),
primitive permutation (\htmladdnormallink{A19}{http://www.research.att.com/cgi-bin/access.cgi/as/njas/sequences/eisA.cgi?Anum=A000019}),
transitive permutation (\htmladdnormallink{A2106}{http://www.research.att.com/cgi-bin/access.cgi/as/njas/sequences/eisA.cgi?Anum=A002106}),
simple (\htmladdnormallink{A5180}{http://www.research.att.com/cgi-bin/access.cgi/as/njas/sequences/eisA.cgi?Anum=A005180}),
total number (\htmladdnormallink{A1}{http://www.research.att.com/cgi-bin/access.cgi/as/njas/sequences/eisA.cgi?Anum=A000001}), $\ldots$.

In the past year the main look-up page has been translated into 28
languages, with the goal of making it easier to use throughout the world.
The entries from the database still appear in English, but the headings in the replies and the error messages have also been translated.

\section{The database}
To be included in the database a sequence should be integer-valued, well-defined and interesting.

The main sources are combinatorics, number theory and recreational mathematics, but most branches of mathematics are represented (e.g. \htmladdnormallink{A27623}{http://www.research.att.com/cgi-bin/access.cgi/as/njas/sequences/eisA.cgi?Anum=A027623}, the number of
rings with $n$ elements), and there are hundreds of entries from chemistry and physics
(e.g. \htmladdnormallink{A8253}{http://www.research.att.com/cgi-bin/access.cgi/as/njas/sequences/eisA.cgi?Anum=A008253}, the coordination sequence for diamond:
the number of carbon atoms that are $n$ bonds away from a particular carbon atom).

Sequences of rational numbers are entered as a linked pair giving numerators and denominators separately.
The Bernoulli numbers $B_n$ are \htmladdnormallink{A27641}{http://www.research.att.com/cgi-bin/access.cgi/as/njas/sequences/eisA.cgi?Anum=A027641}/\htmladdnormallink{A27642}{http://www.research.att.com/cgi-bin/access.cgi/as/njas/sequences/eisA.cgi?Anum=A027642}.

Triangular arrays of numbers are read row-by-row,
so that Pascal's triangle gives \htmladdnormallink{A7318}{http://www.research.att.com/cgi-bin/access.cgi/as/njas/sequences/eisA.cgi?Anum=A007318}:
$$1,~1,1,~1,2,1,~1,3,3,1,~1,4,6,4,1,~ \ldots ~.$$
Square arrays are read by anti-diagonals, so that the Nim-addition table
$$
\begin{array}{ccccccc}
0 & 1 & 2 & 3 & 4 & 5 & \ldots \\
1 & 0 & 3 & 2 & 5 & 4 & \ldots \\
2 & 3 & 0 & 1 & 6 & 7 & \ldots \\
3 & 2 & 1 & 0 & 7 & 6 & \ldots \\
\multicolumn{2}{c}{\ldots} & \multicolumn{2}{c}{\ldots} &
\multicolumn{2}{c}{\ldots} & \ldots
\end{array}
$$
produces \htmladdnormallink{A3987}{http://www.research.att.com/cgi-bin/access.cgi/as/njas/sequences/eisA.cgi?Anum=A003987}:
$$0,~1,1,~2,0,2,~3,3,3,3,~4,2,0,2,4,~ \ldots \,.$$

Most well-defined submissions get accepted, since an open-door policy seems the best.
The amazing coincidences of the Monstrous Moonshine investigations \cite{CN79},
for example, make it difficult to say that a particular sequence, no matter how obscure, will never be of interest.

Sequences that are discouraged are those that depend on an arbitrary and large parameter:
primes of the form $x^2 + y^2 + 2003$, say,
whereas primes of the form $x^2 + y^2 +1$ form a perfectly acceptable sequence
(\htmladdnormallink{A79545}{http://www.research.att.com/cgi-bin/access.cgi/as/njas/sequences/eisA.cgi?Anum=A079545}).

The Encyclopedia currently receives between $10,000$ and $12,000$ downloads per day.
The rate of arrival of {\em new} sequences has remained constant at about $10,000$
per year for the past seven years, with roughly the same number of comments
and updates.
To keep this flood of information from getting out of control, people are
asked to use a 
web form\footnote{\htmladdnormallink{http://www.research.att.com/$\sim$njas/sequences/submit.html}{www.research.att.com/~njas/sequences/submit.html}}
when submitting new sequences or comments.

For most of its life the Encyclopedia has been maintained
by the author, but in the past year a board of
associate editors has been formed to help with the work.
There is also a group of regular users who constantly send corrections and
extensions, and help maintain the accuracy of the entries.

Even so, much remains to be done.
There are more journals and e--print servers now than ever,
and the trained eye sees integer sequences everywhere.
I still discover articles in the library or on the web where authors have
published sequences without sending them to the Encyclopedia.
If you come across an integer sequence in your own work or elsewhere, please 
submit it to the Encyclopedia!

Of course, accuracy is a major concern in maintaining the database.
The entries in \cite{HIS} and \cite{EIS} were checked very thoroughly,
and almost all the errors that have been discovered in those
books were already present in the sources from
which the sequences were taken.
As the number of sequences has increased in
recent years, it has become more difficult to check them all.
However, the number of users has also increased,
and a large number of the entries carry a comment
that the sequence has been extended (or sometimes,
corrected and extended) by someone.
Contributors see a reminder that the standards are those of a mathematics
reference work, and all submissions should be carefully checked.
So, on the whole, users can be confident that the
sequences are correct.
The keywords ``uned'' and ``obsc'' indicate sequences
that have not yet been edited or for
which the definition is unclear.
These serve both as warnings to users and
to indicate places where volunteers can help.

One of the pleasures of maintaining the database
is seeing the endless flow of new sequences.
I will end by mentioning a few
recent examples:  

{\bf Home primes} (\htmladdnormallink{A37274}{http://www.research.att.com/cgi-bin/access.cgi/as/njas/sequences/eisA.cgi?Anum=A037274}), \cite{Gee02}:
$a(n)$ is the prime reached when you start with $n$, concatenate its
prime factors, and repeat until a prime is reached
($a(n)$ is defined to be $-1$ if no prime is ever reached,
although it is conjectured that this never happens).
E.g. $8 = 2 \times 2 \times 2 \rightarrow 222
= 2 \times 3 \times 37 \rightarrow 2337 = 3 \times 19 \times 41 
\rightarrow 31941 \rightarrow \ldots \rightarrow \mbox{~(after ~ 13 ~ steps)}~ 3331113965338635107$,
a prime; so $a(8) = 3331113965338635107$:
$$1~~2~~3~~211~~5~~23~~7~~3331113965338635107~~311~~773~~11~~223~~13~~\ldots .$$

{\bf The EKG sequence} (\htmladdnormallink{A64413}{http://www.research.att.com/cgi-bin/access.cgi/as/njas/sequences/eisA.cgi?Anum=A064413}), \cite{EKG}: 
$a(1) =1$, $a(2) =2$ and, for $n \ge 3$, $a(n)$ is the smallest
natural number not in $\{a(k) : 1 \le k \le n-1\}$ with the
property that ${\rm gcd} \{a(n-1), a(n) \} \ge 2$:
$$1~~2~~4~~6~~3~~9~~12~~8~~10~~5~~15~~18~~14~~7~~21~~24~~16~~20~~22~~11~~\ldots .$$

{\bf Lacing a shoe} (\htmladdnormallink{A78601}{http://www.research.att.com/cgi-bin/access.cgi/as/njas/sequences/eisA.cgi?Anum=A078601}), \cite{Pol02}:
Number of ways to lace a shoe that has n pairs of eyelets.
The lace must follow a Hamiltonian path through the 2n eyelets, and at least one of the neighbors of every eyelet must be on the other side of the shoe.
$$1~~3~~42~~1080~~51840~~3758400~~382838400~~52733721600~~\ldots .$$

{\bf A ``bootstrap'' sequence} (\htmladdnormallink{A79000}{http://www.research.att.com/cgi-bin/access.cgi/as/njas/sequences/eisA.cgi?Anum=A079000}), \cite{CSV03}:
$a(n)$ is taken to be the smallest positive integer greater than $a(n-1)$
which is consistent with the assertion
``$n$ is a member of the sequence if and only if $a(n)$ is odd.''
$$1~~4~~6~~7~~8~~9~~11~~13~~15~~16~~17~~18~~19~~20~~21~~23~~25~~27~~29~~31~~\ldots .$$


\end{document}